\documentclass[11pt]{article}
\usepackage{amsfonts,latexsym,amsmath,amscd,geometry}
\usepackage{amssymb}
\usepackage{latexsym}

\newcommand \nc{\newcommand}
\newtheorem{theorem}{Theorem}[section]
\newtheorem{lemma}[theorem]{Lemma}
\newtheorem{proposition}[theorem]{Proposition}
\newtheorem{corollary}[theorem]{Corollary}

\newtheorem{remark}[theorem]{Remark}

\nc{\ba}{\begin{array}}\nc{\ea}{\end{array}}
\nc{\be}{\begin{eqnarray}}\nc{\ee}{\end{eqnarray}}
\nc{\beq}{\begin{equation}}\nc{\eeq}{\end{equation}}
\nc{\bex}{\begin{eqnarray*}}\nc{\eex}{\end{eqnarray*}}
\nc{\btm}{\begin{theorem}} \nc{\etm}{\end{theorem}}
\nc{\blm}{\begin{lemma}} \nc{\elm}{\end{lemma}}
\nc{\R}{\mathbb{R}} \nc{\va}{\varepsilon} \nc{\ls}{\limits}

\def\pf{\noindent{\bf Proof.\quad}}\def\endpf{\hfill$\Box$}

\newcommand \qed {\hfill $\Box$}

\begin{document}
\title{On uniqueness of heat flow of harmonic maps}

\author{Tao Huang \quad and\quad Changyou Wang\\
\\
Department of Mathematics, University of Kentucky\\
Lexington, KY 40506, USA}

\date{}
\maketitle

\begin{abstract} In this paper, we establish the uniqueness of heat flow of harmonic maps into $(N,h)$ that have sufficiently
small renormalized energy, provided that $N$ is either a unit sphere $S^{k-1}$ or a compact Riemannian homogeneous manifold without boundary.
For such a class of solutions, we also establish the convexity property
of the Dirichlet energy for $t\ge t_0>0$ and the unique limit property at time infinity. As a corollary,
the uniqueness is shown for heat flow of harmonic maps into any compact Riemannian manifold $N$ without boundary
whose gradients
belong to $L^q_t L^l_x$,  for $q>2$ and $l>n$ satisfying (\ref{serrin_condition}).
\end{abstract}

\section {Introduction}
\setcounter{equation}{0}
\setcounter{theorem}{0}

It is well known that for geometric nonlinear evolution equations
with critical nonlinearity, the uniqueness and regularity of weak
solutions is often a very challenging question. In this paper, we
aim to address the issue of uniqueness for heat flow of harmonic maps
in dimensions $n\ge 2$.

Let $(M,g)$ be a $n$-dimensional compact Riemannian manifold
possibly with $\partial M\not=\emptyset$
or complete Riemannian manifold with $\partial M=\emptyset$, and $(N,h)\subset\mathbb R^k$ be a compact
Riemannian manifold without boundary. For $0<T\le +\infty$, the heat flow of harmonic maps for $u:M\times [0,T)\rightarrow N$ is:
\begin{equation}\label{heatflow}
\begin{cases}
\partial_t u-\Delta_g u=A(u)(\nabla u,\nabla u) &\ {\rm{in}}\ M\times (0,T) \\
u=u_0 & \ {\rm{on}}\ \partial_p(M\times [0,T])
\end{cases}
\end{equation}
where $\Delta_g$ is the Laplacian operator on $(M,g)$,
$A(\cdot)(\cdot,\cdot)$ is the second fundamental form of $N\subset\mathbb R^k$,
$\partial_p (M\times [0,T])=(M\times \{0\})\cup(\partial M\times (0,T])$ denotes the parabolic boundary of $M\times [0,T]$,
and $u_0:M\rightarrow N$ is a given map.

The heat flow of harmonic maps has been extensively studied in the past several decades. Under certain geometric conditions on $(N,h)$, the existence of a unique, global smooth solution to (\ref{heatflow})  has been established by Eells-Sampson \cite{ES},
Hamilton \cite{H}, and Hildebrandt-Kaul-Widman \cite{HKW}. In general, the existence of a unique, global weak
solution to (\ref{heatflow}) with finitely many singularities has been obtained by Struwe \cite{struwe} and
Chang \cite{chang}  for $n=2$; and the existence of partially regular, global weak
solutions to (\ref{heatflow}) has been established by Chen-Struwe \cite {chen-struwe} and Chen-Lin \cite{chen-lin} for $n\ge 3$.
Concerning the uniqueness for weak solutions to (\ref{heatflow}),  Freire \cite{freire} first proved that in dimension $n=2$,
the uniqueness holds for weak solutions whose Dirichlet energy is monotone decreasing with respect to $t$
(see  L.Wang \cite{luwang} and L. Z. Lin \cite{lin} for a new simple proof). For $n\ge 3$, there are non-uniqueness for weak solutions to (\ref{heatflow}), see the examples constructed by Coron \cite{coron} and Bethuel-Coron-Ghidaglia-Soyeur \cite{BBC}. In fact, Coron \cite{coron} proved that for suitable initial data, there exist weak solutions to (\ref{heatflow}) that are different from those constructed
by  Chen-Struwe. Partially motivated
by \cite{coron}, Struwe \cite{struwe3} has raised the following question:\\
 {\it For $M=\mathbb R^n$, exhibit  a class of functions within which (\ref{heatflow}) posses a unique solution. Certainly the class of functions satisfying the strong monotonicity
formula
\begin{equation}
\label{struwe_mono0}
\Phi_{(\bar x, \bar t)}(\rho)\le\Phi_{(\bar x, \bar t)}(r),  \ \forall \bar x\in\mathbb R^n, \bar t>0,
\ 0<\rho\le r\le\sqrt{\bar t}
\end{equation}
is a likely candidate. Here
$$\Phi_{(\bar x,\bar t)}(\rho)=\rho^2\int_{\mathbb R^n\times\{\bar t-\rho^2\}}|\nabla u|^2(x,t) G(x-\bar x, t-\bar t)\,dx$$
and
$$G(y,s)=\frac{1}{(4\pi|s|)^{\frac{n}2}}\exp\Big(-\frac{|y|^2}{4|s|}\Big), \ y\in \mathbb R^n, s<0 $$
is the fundamental solution to the backward heat equation on $\mathbb R^n$.}\\
To the best of authors' knowledge, this question is largely open.
In this paper, we will obtain  some uniqueness results
for the heat flow of harmonic maps (\ref{heatflow}), that may  shed
light on the validity of Struwe's conjecture as above.

To state the result, we introduce some notations. For $1<p, q<\infty$
and $0<T\le +\infty$, define the Sobolev
spaces
$$W^{1,p}(M,N):=\left\{v\in W^{1,p}(M,\R^k)\ \Big | \ v(x)\in N, \mbox{ a.e. } x\in M\right\},$$
$$H^{1}\Big(M\times[0,T],N\Big):=\left\{v\in W^{1,2}(M\times[0,T],\R^k)
\ \Big|\ v(x,t)\in N, \mbox{ a.e. } (x,t)\in M\times[0,T]\right\},$$
the $L^q_tL^p_x$-space
$$L^q_t L^p_x\Big(M\times[0,T], \mathbb R^k\Big)
:=\Big\{ f:M\times [0,T]\to\mathbb R^k
\ \Big| \ f\in L^q([0,T], L^p(M))\Big\},$$
and the parabolic Morrey space ${M}^{p, \lambda}_{R}$.  For any $1\leq p<+\infty$,
$0\leq\lambda\leq n+2$,  $0<R\le+\infty$, and any open set $U=U_1\times U_2\subset M\times\mathbb R$,
$${M}^{p,\lambda}_R(U)=\left\{f\in L^{p}_{\mbox{loc}}(U):
\Big\|f\Big\|_{{M}^{p,\lambda}_R(U)}<+\infty\right\},$$
where
$$\Big\|f\Big\|_{{M}^{p,\lambda}_R(U)}=\Big(\sup\limits_{(x,t)\in U}\sup\limits_{0<r<\min\{R, d_g(x,\partial U_1),\sqrt{t}\}}\
r^{\lambda-n-2}\int_{P_r(x,t)\cap U} |f|^p\Big)^{\frac{1}{p}}. $$
Here $d_g$ denotes the distance function on $M$ induced by $g$, and
$$ P_r(x,t)=B_r(x)\times (t-r^2,t], \ {\rm{with}}\ B_r(x)=\{y\in M: \ d_g(y,x)\le r\}$$
for $(x,t)\in U$ and $0<r<d_g(x,\partial U_1):=\displaystyle\inf_{y\in\partial U_1} d_g(x,y)\footnote{if $\partial U=\emptyset$,
then we set $d_g(x,\partial U)=\infty$.}$. When $R=\infty$, we denote $M^{p,\lambda}(U)=M^{p,\lambda}_\infty(U)$.

For $u_0\in W^{1,2}(M,N)$ and $0<T\leq +\infty$,
$u\in H^1(M\times[0,T],N)$ is a weak solution of
(\ref{heatflow}) if $u$ satisfies
$(\ref{heatflow})_1$ in the sense of distribution and
$(\ref{heatflow})_2$ in the sense of trace.

Now we state our main theorem on the uniqueness of weak solutions to (\ref{heatflow}).
\begin{theorem} \label{unique1} For $n\ge 2$ and $1<p\le 2$, there exist $\epsilon_0=\epsilon_0(p,n)>0$ and $R_0=R_0(M, g,\epsilon_0)>0$ such that if \\
(i) $(M,g)$ is a $n$-dimensional Riemannian manifold that is either complete noncompact without boundary or compact with or without boundary;\\
(ii) $(N,h)\subset\mathbb R^k$ is either the unit sphere $S^{k-1}$ or a compact Riemannian
homogeneous manifold without boundary; and \\
(iii) $u_1,u_2\in H^1(M\times [0,T],N)$ are two weak
solutions of (\ref{heatflow}), with $u_1=u_2=u_0$ on $\partial_p (M\times [0,T])$ for some $u_0\in W^{1,2}(M,N)$,  that satisfy
\begin{equation}\label{struwe_mono}
\max_{i=1,2} \Big[\|\nabla u_i\|_{M^{p,p}_{R_0}(M\times (0,T))}+\|\partial_t u_i\|_{M^{p,2p}_{R_0}(M\times (0,T))}\Big]
\leq \epsilon_0,
\end{equation}
then $u_1\equiv u_2$ on $M\times [0,T]$.
\end{theorem}

Recall that $N$ is a Riemannian homogeneous manifold if there exists a finite dimensional Lie group $\mathcal G$ (dim $\mathcal G=s<+\infty$) that acts
transitively on $N$ by isometries.

There are two main ideas of proof of Theorem \ref{unique1}: \\
(i) an $\epsilon_0$-regularity theorem (Theorem \ref{e-regularity} in \S2 below) for the heat flow
of harmonic maps that satisfy the smallness condition (\ref{struwe_mono}), which is new and improves
the regularity theorem previously obtained by Chen-Li-Lin \cite{CLL}, Feldman \cite{feldman}, and Chen-Wang \cite{CW}.
It shall have its own interest. In particular, we have that  for $i=1,2$,
$u_i\in C^\infty(M\times (0, T])$ and  satisfies the gradient
estimate:
\begin{equation}
\max_{i=1,2}|\nabla u_i|(x,t)
\le C\epsilon_0\Big(\frac{1}{R_0}+\frac{1}{d_g(x,\partial M)}+\frac{1}{\sqrt{t}}\Big),
\ \forall (x,t)\in M\times (0,T], \label{gradient_est11}
\end{equation}
(ii) applications of (\ref{gradient_est11}), the Hardy inequality,  and a generalized Gronwall inequality type argument.

Now a few remarks are in order.
\begin{remark}{\rm
i) Due to technical difficulties,  it is unknown whether the $\epsilon$-regularity Theorem \ref{e-regularity} (with $p=2$) holds for a general Riemannian manifold $N$.
Hence it is an open question that Theorem \ref{unique1}, Theorem \ref{convexity1}, and Corollary \ref{unique_limit} remain to hold
for a general manifold $N$. \\
ii) Note that by the H\"older inequality, the Morrey norm  $\mathcal E(p):=(\displaystyle\|\nabla u\|_{M^{p,p}(\cdot)}+\|\partial_t u\|_{M^{p,2p}(\cdot)})$
is monotone increasing for $1<p\le 2$. The bound  of $\mathcal E(2)$
for solutions $u$ to (\ref{heatflow})  holds if $u$ satisfies\\
(a) a local energy inequality  (assume
$M=\mathbb R^n$ for simplicity):
\begin{equation}\label{local_energy_ineq}
\int_{P_{r}(x,t)}|\partial_t u|^2
\le\frac{C}{(R-r)^2} \int_{P_R(x,t)}|\nabla u|^2, \quad \forall (x,t)\in \mathbb R^{n+1}_+,
 \ 0<r\le\frac{R}2, \ R\le\sqrt{t},
\end{equation}
(b) a local energy monotonicity inequality:
\begin{equation}\label{local_struwe}
r^{-n}\int_{P_{r}(x,t)}|\nabla u|^2 \le C R^{-n}\int_{P_R(x,t)}|\nabla u|^2, \quad \forall (x,t)\in \mathbb R^{n+1}_+,
\ \ 0<r\le \frac{R}2, \ R\le\sqrt{t}.
\end{equation}
Both properties hold if $u$ is either a smooth solution (see \cite{struwe} and \cite{chen-struwe}) or a stationary solution
of (\ref{heatflow}) (see \cite{CLL}, \cite{feldman}, and \cite{CW}).
Therefore, under (\ref{local_energy_ineq}) and (\ref{local_struwe}),
the condition (\ref{struwe_mono}) is satisfied, provided that there exists $R_0>0$ such that
there holds
\begin{equation}\label{struwe_mono1}
\sup\Big\{
R_1^{-n}\int_{P_{R_1}(x,t)}|\nabla u|^2
\ \Big| \ x\in\mathbb R^n, \ R_1=\min\{R_0, \sqrt{t}\}\Big\}\le \epsilon^2_0.
\end{equation}
Hence Theorem \ref{unique1} implies that the uniqueness does hold for the class of solutions that satisfy, in addition
to (\ref{local_energy_ineq})  and (\ref{local_struwe}), the smallness condition (\ref{struwe_mono1}).\\
iii)  For any compact or complete noncompact $(M, g)$ without boundary, there exists $\epsilon_0>0$
such that if the initial data $u_0:M\to N$ satisfies that for some  $R_0>0$,
\begin{equation}\label{small_initial}
\sup\Big\{ r^{2-n}\int_{B_r(x)}|\nabla u_0|^2
\ \Big| x\in M, \ r\le R_0\Big\} \le \epsilon_0^2,
\end{equation}
then as a consequence of the local well-posedness theorem by Wang \cite{wang-well}, there exists $0<T_0(\approx R_0^2)$ and
a solution $u\in C^\infty(M\times (0,T_0),N)$ of (\ref{heatflow}) that satisfies the condition (\ref{struwe_mono}).}
\end{remark}

Motivated by the proof of Theorem \ref{unique1}, we obtain the following convexity property on (\ref{heatflow}).
\begin{theorem} \label{convexity1}
For $n\ge 2$, $1<p\le 2$, and $1\le T\le\infty$,  there exist $\epsilon_0=\epsilon_0(p,n)>0$, $R_0=R_0(M, g,\epsilon_0)>0$, and
$0<T_0=T_0(\epsilon_0)<T$ such that if \\
(i) $(M,g)$ is a $n$-dimensional Riemannian manifold that is either complete noncompact without boundary or compact with or without boundary;\\
(ii) $(N,h)\subset\mathbb R^k$ is either the unit sphere $S^{k-1}$ or a compact Riemannian
homogeneous manifold without boundary; and \\
(iii) $u\in H^1(M\times [0,T],N)$ is a weak solution of (\ref{heatflow}), with $u=u_0$ on $\partial_p (M\times [0,T])$ for some $u_0\in W^{1,2}(M,N)$,  that satisfies
\begin{equation}\label{struwe_mono00}
\|\nabla u\|_{M^{p,p}_{R_0}(M\times (0,T))}+\|\partial_t u\|_{M^{p,2p}_{R_0}(M\times (0,T))}
\leq \epsilon_0,
\end{equation}
then \\
(i) the Dirichlet energy $E(u(t)):=\frac12\int_M |\nabla u|^2$ is monotone decreasing for $t\ge T_0$;
and \\
(ii) for any $t_2\ge t_1\ge T_0$,
\begin{equation}\label{convexity2}
\int_{M}|\nabla ( u(t_1)-u(t_2))|^2\le C\Big[\int_{M}|\nabla u(t_1)|^2-\int_M |\nabla u(t_2)|^2\Big].
\end{equation}
\end{theorem}

We would like to remark that the convexity property has been observed by Schoen \cite{schoen}
for the Dirichlet energy of harmonic maps into manifolds $N$ with nonpositive sectional curvatures. In \S5 appendix below,
we will show that it also holds
for harmonic maps with small renormalized energy, which yields a new proof of the uniqueness theorem
by Struwe \cite{struwe2} and Moser \cite{moser1}.

A direct consequence of Theorem \ref{convexity1} is the following uniqueness of limit at $t=\infty$ for (\ref{heatflow}).
\begin{corollary}\label{unique_limit}
For $n\ge 2$ and $1<p\le 2$,  there exist $\epsilon_0=\epsilon_0(p,n)>0$, and $R_0=R_0(M, g,\epsilon_0)>0$ such that if \\
(i) $(M,g)$ is a $n$-dimensional Riemannian manifold that is either complete noncompact without boundary or compact with or without boundary;\\
(ii) $(N,h)\subset\mathbb R^k$ is either the unit sphere $S^{k-1}$ or a compact Riemannian
homogeneous manifold without boundary; and \\
(iii) $u\in H^1(M\times [0,\infty),N)$ is a weak solution of (\ref{heatflow}), with $u=u_0$ on $\partial_p (M\times [0,\infty])$ for some $u_0\in W^{1,2}(M,N)$,  that satisfies the condition (\ref{struwe_mono00}).\\
Then there exists a harmonic map $u_\infty\in C^\infty(M, N)\cap W^{1,2}(M, N)$, with $u_\infty=u_0$ on $\partial M$,
such that
\begin{equation}\label{unique_limit1}
\lim_{t\uparrow\infty}\|u(t)-u_\infty\|_{W^{1,2}(M)}=0,
\end{equation}
and, for any compact subset $K\subset\subset M$ and $m\ge 1$,
\begin{equation}\label{unique_limit2}
\lim_{t\uparrow\infty}\|u(t)-u_\infty\|_{C^m(K)}=0.
\end{equation}
\end{corollary}

The uniqueness of limit at $t=\infty$ has been proved by Hartman \cite{hartman} for the smooth solutions to (\ref{heatflow}) when $N$ has
nonpositive sectional curvatures. L. Simon in his celebrated work \cite{simon} has shown the unique limit at $t=\infty$
for smooth solutions to (\ref{heatflow}) into a target manifold $(N,h)$ that is real analytic. Note that the solution $u$ in Theorem \ref{convexity1}
is allowed to be singular near the parabolic boundary $\partial_p(M\times (0,\infty))$, as the initial-boundary data $u_0$ may be
in $W^{1,2}(M,N)$.  Also, our proof of Theorem \ref{convexity1} depends only on the smallness condition (\ref{struwe_mono00}) and the small energy
regularity Theorem \ref{e-regularity}.
During the preparation of this work, we have seen two very interesting articles by L.Wang \cite{luwang} and L.Z. Lin \cite{lin}, in which
 Theorems \ref{unique1}, \ref{convexity1}, and Corollary \ref{unique_limit}  were proven for Struwe's almost regular solution $u$ to (\ref{heatflow})
in dimension $n=2$  when the Dirichlet energy of $u_0$ is sufficiently small. We would like to point that
since Struwe's solution $u$ to (\ref{heatflow}) satisfies the energy inequality, the condition in \cite{lin} yields the global smallness:
$$\sup_{t\ge 0} E(u(t))+\int_{M\times [0, t]}|\partial_t u|^2\le E(u_0)\le\epsilon_0^2, \ \forall t>0,$$
which is stronger than (\ref{struwe_mono00}) in dimension $n=2$.
There is also an interesting paper by Topping \cite{topping}
that addressed the rigidity at $t=\infty$ of the heat flow of harmonic maps from $S^2$ to $S^2$.

A class of weak solutions that satisfy the smallness condition (\ref{struwe_mono00}) are the
so-called Serrin's $(l,q)$-solutions.
Recall that a weak solution $u\in H^1_{\hbox{loc}}(M\times [0,T],N)$ of (\ref{heatflow}) is called a Serrin's $(l,q)$-solution if, in additions,
$\nabla u\in L^{q}_t L^{l}_x(M\times [0,T])$ for some $l\ge n$ and $q\geq 2$ satisfying
\beq\label{serrin_condition}
\frac{n}{l}+\frac{2}{q}=1. \eeq

In \S3 below, we will verify that if $u$ is a Serrin's $(l,q)$-solution to (\ref{heatflow}) with $l>n$, and $\displaystyle u|_{\partial_p (M\times [0,T])}=u_0$ for
a given $u_0:M\to N$ with $\nabla u_0\in L^r(M)$ for some $n<r<\infty$,
then $u$ satisfies (\ref{struwe_mono00}) for some $p_0>1$. We would also like to point
out for such an initial and boundary data $u_0$, the local existence of Serrin's $(l,q)$-solutions of (\ref{heatflow}) can be shown by the standard fixed point
theory. In fact, interested readers may verify that the argument by Fabes-Jones-Riviere \cite{FJR} \S4 can be adapted to achieve such an existence.
Here we have the following uniqueness result for Serrin's ($p,q$)-solutions of the heat flow of harmonic maps into a general Riemannian manifold.

\begin{theorem}\label{unique2}
For $n\geq 2$, $0<T\leq +\infty$, let
$(M,g)$ be either a compact or complete Riemannian manifold without boundary or a compact Riemannian manifold with boundary,
and $N$ be a  compact Riemannian manifold without boundary.
Let $u_1,u_2\in H^1(M\times [0,T],N)$ be two weak solutions of (\ref{heatflow}),
with $u_1=u_2=u_0$ on $\partial_p(M\times [0,T])$ for some $u_0\in W^{1,2}(M,N)$, such that $\nabla u_1,\nabla u_2\in L^q_tL^l_x(M\times [0,T])$
for some $(l,q)$ satisfying (\ref{serrin_condition}) with $l>n, q>2$.
Then $u_1, u_2\in C^\infty(M\times (0, T))$, and $u_1\equiv u_2$ on $M\times[0,T]$.
\end{theorem}

\begin{remark}{\rm
We conjecture that Theorem \ref{unique2} remains to be true for the end point case
$l=n, q=+\infty$. We would like to point out that Lin-Wang \cite{lin-wang}
have proved the uniqueness holds for two weak solutions $u_1, u_2$ to (\ref{heatflow})
with the same initial data, provided that
$\nabla u_1, \nabla u_2\in C([0, T], L^n(M))$.
Wang \cite{wang} has proved that for any $n\geq 4$, any weak solution $u\in H^1(M\times [0,T], N)$
with $\nabla u\in L^{\infty}_t L^n_x(M\times [0,T])$ belongs to $C^\infty(M\times (0,T])$.
However, since
$\|\nabla u(t)\|_{L^n(M)}$ may lack continuity at  $t=0$,
the issue of uniqueness for the end point case remains unsolved.}
\end{remark}

It turns out that we can extend the ideas in this paper to study the uniqueness issue of heat flow of biharmonic maps,
which will be discussed in a forthcoming paper \cite{HHW}.

The paper is written as follows. In \S2, we will provide an $\epsilon$-regularity theorem on
certain weak solutions of (\ref{heatflow}) for $N$ either a unit sphere or a compact Riemannian homogeneous manifold
without boundary.
In \S3, we will outline the proofs of Theorem \ref{unique1}, Theorem \ref{convexity1}, and
Corollary \ref{unique_limit}.  In \S4, we will discuss Serrin's $(l,q)$-solutions of (\ref{heatflow}) and
sketch a proof of Theorem \ref{unique2}.  In \S5, we will provide
a simple alternative proof to an improved version of the uniqueness theorem for harmonic maps with small energy,
originally due to Struwe \cite{struwe2} ($n=3$) and Moser \cite{moser1} ($n\ge 4$).

\section{$\epsilon$-regularity Theorem}
\setcounter{equation}{0}
\setcounter{theorem}{0}

In this section, we will establish an $\epsilon$-regularity theorem for the heat flow of harmonic maps (\ref{heatflow}),
which plays a crucial role in the proof of our main theorems. This regularity theorem seems to be new, whose proof is rather
elementary.  It improves the regularity theorem previously obtained by Chen-Li-Lin \cite{CLL}, Feldman \cite{feldman},
Chen-Wang \cite{CW} (see also Moser \cite{moser, moser3} for more general results).  We believe
that it shall have its own interests. We would also like to point out the relevant works on the regularity theorem
on stationary harmonic maps by H\'elein \cite{helein0}, Evans \cite{evans}, Bethuel \cite{bethuel}, Chang-Wang-Yang
\cite{CWY}, and Riviere-Struwe \cite{RS}. Especially, the proof of the regularity theorem \ref{e-regularity}
below is motivated by  \cite{CWY}.

From now on, let $\Omega\subset\mathbb R^n$ be a bounded smooth domain, and denote
$$\delta\Big((x,t), (y,s)\Big)=\max\Big\{|x-y|, \ \sqrt{|t-s|}\Big\},
\ (x,t), \ (y,s)\in \mathbb R^n\times \mathbb R$$
as  the parabolic distance on $\mathbb R^n\times \mathbb R$.

\begin{theorem} \label{e-regularity} Assume that $N$  is either a unit sphere $S^{k-1}$ or a compact Riemannian homogeneous manifold without boundary.
For $1<p\le 2$ and $0<T<+\infty$, there exists $\epsilon_p>0$  such that if
$u\in H^1(\Omega\times [0,T], N)$ is a weak solution of (\ref{heatflow})$_1$ and  satisfies that,
for $z_0=(x_0, t_0)\in\Omega\times (0, T]$ and  $0<R_0\le \frac12 \min\{d(x_0,\partial\Omega), \sqrt{t_0}\}$,
\begin{equation}\label{small_norm1}
\|\nabla u\|_{M^{p,p}_{R_0}(P_{R_0}(z_0))}+\|\partial_t u\|_{M^{p,2p}_{R_0}(P_{R_0}(z_0))}
\le\epsilon_p.
\end{equation}
Then $u\in C^\infty(P_{\frac{R_0}4}(z_0), N)$, and
\begin{equation}\label{gradient_est1}
|\nabla^m u|(z_0)\le \frac{C\epsilon_p}{R_0^m}, \ \forall m\ge 1.
\end{equation}
\end{theorem}
\begin{remark} It remains an open question whether Theorem \ref{e-regularity} holds for  any
compact Riemannian manifold $N$ without boundary, under the condition (\ref{small_norm1}) for $p=2$. The interested
readers can refer to Moser \cite{moser2} and Moser \cite{moser3} for related works.
\end{remark}

The proof of Theorem \ref{e-regularity} is based on the following lemma.

\begin{lemma} \label{grad_estimate5}
For any $1<p\le 2$, there exists $\epsilon_p>0$ such that if $N=S^{k-1}$ or a compact Riemannian
homogeneous manifold without boundary, and
$u\in H^1(P_4, N)$  is a weak solution of (\ref{heatflow}) satisfying
\begin{equation}
\sup_{(x,t)\in P_2, 0<r\le 2} r^{p-(n+2)}\int_{P_r(x,t)}(|\nabla u|^p+r^p|\partial_tu|^p)
\le\epsilon^p.\label{small_normal}
\end{equation}
Then $u\in  C^\infty(P_{\frac12}, S^{k-1})$ and satisfies
\begin{equation}\label{grad_estimate6}
\|\nabla^m u\|_{C^0(P_{\frac12})}\le C(n, p, \epsilon, m), \ \forall m\ge 1.
\end{equation}
\end{lemma}
\pf The crucial step to establish (\ref{grad_estimate6}) is the following decay
estimate:\\
\noindent{\bf Claim}: There exists $q>\max\{\frac{p}{p-1}, n+2\}$ such that 
for any $\theta\in (0,\frac12)$, 
$z_0=(x_0,t_0)\in P_1$, and $0<r\le1$, it holds
\begin{equation}\label{decay_estimate1}
\frac{1}{(\theta r)^{n+2}}\int_{P_{\theta r}(z_0)}|u-u_{z_0, \theta r}|
\le C\Big(\theta^{-(n+2)}\epsilon_p+\theta\Big) \Big(\frac{1}{r^{n+2}}\int_{P_{r}(z_0)}|u-u_{z_0, r}|^{q}\Big)^{\frac1{q}},
\end{equation}
where $f_{z_0,r}=\frac{1}{|P_r(z_0)|}\int_{P_r(z_0)} f$ is the average of $f$ over $P_r(z_0)$.

For $z_0=(x_0,t_0)\in P_1$ and $0<r\le 1$, since
$v(y,s)=u(z_0+(ry, r^2s)): P_2\to N$ satisfies (\ref{heatflow}),  and the condition
(\ref{small_normal}) yields that $v$ satisfies
\begin{equation}
\sup_{(x,t)\in P_1, 0<r\le 1} r^{p-(n+2)}\int_{P_r(x,t)}|\nabla v|^p
+r^p|\partial_t v|^p\le\epsilon^p.
\label{small_normal1}
\end{equation}
Thus it suffices to show (\ref{decay_estimate1}) for $z_0=(0,0)$ and $r=2$.

We divide the proof into two cases:\\
\noindent {\bf Case 1}: $N=S^{k-1}$ is the unit sphere.\\
\noindent{\it Step 1}. Rewriting of (\ref{heatflow}). Since $|u|=1$, we have
$u^i u^i_\alpha=0$. Also, it follows (\ref{heatflow}) that
$$(u^iu^j_\alpha-u^ju^i_\alpha)_\alpha
=u^i \Delta u^j-u^j\Delta u^i=u^i \partial_t u^j-u^j \partial_t u^i.$$
Hence we have
\begin{eqnarray}
\partial_t u^i-\Delta u^i&=& A^i(u)(\nabla u, \nabla u)=|\nabla u|^2 u^i\nonumber\\
&=&u^j_\alpha u^j_\alpha u^i-u^j_\alpha u^j u^i_\alpha
=u^j_\alpha(u^i u^j_\alpha-u^j u^i_\alpha)\nonumber\\
&=& \Big[(u^j-c^j)(u^iu^j_\alpha-u^ju^i_\alpha)\Big]_\alpha
-(u^j-c^j)(u^i \partial_t u^j-u^j\partial_t u^i), \label{hfhms1}
\end{eqnarray}
where $c^j\in\mathbb R$ is an arbitrary constant. For the convenience, set
$$W^{ij}=u^i\partial_t u^j-u^j\partial_t u^i,\
V^{ij}_\alpha=u^iu^j_\alpha-u^ju^i_\alpha,
\ 1\le i, j\le k, 1\le \alpha\le n.$$
\noindent{\it Step 2}. Construction of auxiliary functions. Let $\eta\in C^\infty_0(P_2)$
such that
$$0\le \eta\le 1, \ \eta=1\ {\rm{on}}\ P_1, \ {\rm{and}} \ |\nabla\eta|\le C.$$
Define $v, w:\mathbb R^n\times \mathbb R_+\to\mathbb R^{k}$ by
\begin{equation}
\partial_t v^i-\Delta v^i=\Big[\eta^2(u^j-c^j)V^{ij}_\alpha\Big]_\alpha;\  v^i \Big|_{t=0}=0,
\end{equation}\label{v-eqn}
and
\begin{equation}
\partial_t w^i-\Delta w^i=-\eta^2(u^j-c^j)W^{ij}; \ w^i\Big|_{t=0}=0.\label{w-eqn}
\end{equation}
Set $h=u-(v+w): P_1\to\mathbb R^{k}$. Then
\begin{equation}
\partial_t h-\Delta h=0\  \ {\rm{in}}\ \ P_1. \label{h-eqn}
\end{equation}
\noindent{\it Step 3}. Estimation of $v, \ w$, and $u$. By the Duhamel formula, we have that
\begin{eqnarray*}
v^i(x,t)&=&\int_0^t\int_{\mathbb R^n}
H(x-y,t-s) \Big[\eta^2(u^j-c^j)V^{ij}_\alpha\Big]_\alpha(y,s)\\
&=&\int_0^t\int_{\mathbb R^n}
\nabla_xH(x-y,t-s) (\eta^2(u^j-c^j)V^{ij}_\alpha)(y,s),
\end{eqnarray*}
where $H$ denotes the heat kernel on $\mathbb R^n$. Then, as in \cite{huang-wang},
we have
$$|\nabla_x H|(x-y,t-s)\lesssim \delta((x,t), (y,s))^{-(n+1)}, \ \  (x,t), \ (y,s)\in \mathbb R^{n+1}, $$
where $\delta((x,t), (y,s))$ is the parabolic distance on $\mathbb R^{n+1}$.
Hence
$$|v^i|(x,t)\lesssim I_1(\eta^2 |u^j-c^j||V^{ij}_\alpha|)(x,t),$$
where
$$I_1(f)(x,t):= \int_{\mathbb R^{n+1}} \frac{f(y,s)}
{\delta((x,t), (y,s))^{n+1}}, \ \forall f\in L^1_{\rm{loc}}(\mathbb R^{n+1}), $$
is the parabolic Riesz potential of order $1$. By the Riesz potential estimate (see \cite{huang-wang}), we have
\begin{eqnarray}\label{v-estimate}
&&\|v\|_{L^1(P_2)}\lesssim\|v\|_{L^p(P_2)}\le\|v\|_{L^{p}(\mathbb R^{n+1})}
\nonumber\\
&\lesssim& \sum_{i,j}\|\eta^2|u^j-c^j||V^{ij}_\alpha|\|_{L^{\frac{(n+2)p}{n+2+p}}(\mathbb R^{n+1})}
\lesssim \sum_{i,j}\|V^{ij}_\alpha\|_{L^p(P_2)}\|u^j-c^j\|_{L^{n+2}(P_2)}.
\end{eqnarray}

For $w$, since
$$w^i(x,t)=\sum_{j}\int_0^t\int_{\mathbb R^n}H(x-y,t-s)(\eta^2(u^j-c^j)W^{ij})(y,s),$$
applying the Young inequality we  obtain
\begin{eqnarray}\label{w-estimate}
\|w\|_{L^1(P_2)}&\le&\|w\|_{L^{\widetilde {q_1}}(P_2)}\le\|w\|_{L^{\widetilde {q_1}}(\mathbb R^n\times [0,1])}\nonumber\\
&\lesssim& \sum_{i,j}\Big\|\eta^2(u^j-c^j)W^{ij}\Big\|_{L^{\widetilde {q_1}}(\mathbb R^n\times [0,1])}
\lesssim \sum_{i,j}\Big\||u^j-c^j||W^{ij}|\Big\|_{L^{\widetilde {q_1}}(P_2)}\nonumber\\
&\lesssim&\sum_{i,j} \Big\|u^j-c^j\Big\|_{L^{q_1}(P_2)}\Big\|W^{ij}\Big\|_{L^p(P_2)},
\end{eqnarray}
where $1<\widetilde {q_1}<p$ and  $q_1>\frac{p}{p-1}$ satisfy
$\displaystyle\frac{1}{\widetilde {q_1}}=\frac{1}{p}+\frac{1}{q_1}.$

For $h$, by the standard theory on the heat equation we have
that for any $0<\theta<1$,
\begin{equation}\label{h-estimate}
\frac1{\theta^{n+2}}\int_{P_\theta}|h-h_\theta|
\lesssim \theta \int_{P_1}|h-h_1|
\lesssim\theta \Big[\|v\|_{L^1(P_2)}+\|w\|_{L^1(P_2)}+\|u-u_2\|_{L^1(P_2)}\Big],
\end{equation}
where $f_r=\frac{1}{|P_r|}\int_{P_r} f$ is the average of a function $f$ over $P_r$.

Now we let $c^j=u^j_2$, the average of $u^j$ over $P_2$ and set 
$q=\max\{q_1, n+2\}$. 
Combining the
estimates (\ref{v-estimate}), (\ref{w-estimate}), and (\ref{h-estimate}) and
applying H\"older's inequality together
yields
\begin{eqnarray}\label{u-estimate}
&&\frac1{\theta^{n+2}}\int_{P_\theta}|u-u_\theta|
\le\frac{1}{\theta^{n+2}}\int_{P_\theta}(|v|+|w|)+\frac1{\theta^{n+2}}\int_{P_\theta}|h-h_\theta|\nonumber\\
&&\lesssim \theta^{-(n+2)} \Big[\|v\|_{L^1(P_2)}+\|w\|_{L^1(P_2)}\Big]+\theta\Big[\|v\|_{L^1(P_2)}+\|w\|_{L^1(P_2)}+\|u-u_2\|_{L^1(P_2)}\Big]\nonumber\\
&&\lesssim \Big[\theta+\theta^{-(n+2)}(\|V^{ij}_\alpha\|_{L^p(P_2)}+\|W^{ij}\|_{L^p(P_2)})\Big]\Big(\frac1{2^{n+2}}\int_{P_2}|u-u_2|^{q}\Big)^{\frac{1}{q}}
\nonumber\\
&&\leq C\Big[\theta+\theta^{-(n+2)}\epsilon_p\Big]\Big(\frac1{2^{n+2}}\int_{P_2}|u-u_2|^{q}\Big)^{\frac{1}{q}},
\end{eqnarray}
where  we have used in the last step
the condition (\ref{small_normal}) so that 
$$\|V^{ij}_\alpha\|_{L^p(P_2)}+\|W^{ij}\|_{L^p(P_2)}\le C\epsilon_p.$$
This yields (\ref{decay_estimate1}). 
It follows from (\ref{small_normal}) and the Poincar\'e
inequality that $u\in {\rm{BMO}}(P_2)$, and
\begin{equation}\label{bmo}\Big[u\Big]_{\rm{BMO}(P_2)}
:=\sup_{P_r(z)\subset P_2}
\Big\{ \frac1{r^{n+2}}\int_{P_r(z)}|u-u_{z,r}|\Big\}\le C\epsilon_p.
\end{equation}
By the celebrated John-Nirenberg's inequality \cite{JN}, (\ref{bmo})  implies that
for any $q>1$, it holds
\begin{equation}\label{reverse}
\sup_{P_r(z)\subset P_2}\Big\{\Big(\frac1{\theta^{n+2}}\int_{P_r(z)}|u-u_{z,r}|^{q}\Big)^{\frac{1}{q}}\Big\}\le C(q)\Big[u\Big]_{\rm{BMO}(P_2)}.
\end{equation}
By (\ref{reverse}), we see that (\ref{decay_estimate1}) implies that
\begin{equation}\label{reverse1}
\frac{1}{(\theta r)^{n+2}}\int_{P_{\theta r}(z_0)}|u-u_{z_0, \theta r}|
\le C\left(\theta^{-(n+2)}\epsilon_p+\theta\right) \Big[u\Big]_{\rm{BMO}(P_2)}
\end{equation}
holds for any $\theta\in (0,\frac12), z_0\in P_1, 0<r\le 1$. Taking supremum of (\ref{reverse1}) over all $z_0\in P_\theta$ and $0<r\le 1$, we obtain
\begin{equation}\label{reverse2}
\Big[u\Big]_{\rm{BMO}(P_\theta)}
\le C\left(\theta^{-(n+2)}\epsilon_p+\theta\right) \Big[u\Big]_{\rm{BMO}(P_2)}.
\end{equation}
If we choose $\theta=\theta_0\in (0,\frac12)$ and $\epsilon_p>0$ so small
that 
$$ C\left(\theta_0^{-(n+2)}\epsilon_p+\theta_0\right)\le\frac12,$$
then (\ref{reverse2}) implies
\begin{equation}\label{reverse3}
\Big[u\Big]_{\rm{BMO}(P_{\theta_0})}
\le \frac12 \Big[u\Big]_{\rm{BMO}(P_2)}.
\end{equation}
It is standard that by iterations and the Campanato theory \cite{cam},  (\ref{reverse3}) implies that there exists $\alpha\in (0,1)$ such that
$u\in C^\alpha(P_{\frac34})$ and
$$\Big[u\Big]_{C^\alpha(P_{\frac34})}\le C(p,\epsilon_p).$$
The higher regularity and the estimate (\ref{grad_estimate6}) then follow from the parabolic hole filling type argument
and the bootstrap argument (see also \cite{huang-wang}).

\medskip
\noindent{\bf Case 2}: $N$ is a compact Riemannian homogeneous manifold without boundary. We will indicate that (\ref{heatflow}) can be written
into the same form as (\ref{hfhms1}). In fact, according to H\'elein \cite{helein}, there exist $s$ smooth tangential vector
fields $Y_1,\cdots, Y_s$ and $s$ smooth tangential killing vector fields $X_1,\cdots, X_s$ on $N$ such that
for any $y\in N$ and $V\in T_y N$, it holds
$$V=\sum_{i=1}^s \langle V, X_i(y)\rangle Y_i(y).$$
Thus, as in \cite {CW} Lemma 4.2, (\ref{heatflow}) is equivalent to
\begin{eqnarray}\label{hfhms2}
\partial_t u-\Delta u &=& -\sum_{i=1}^s \langle \nabla u, X_i(u)\rangle \nabla(Y_i(u))\nonumber\\
&=& -\sum_{i=1}^s \hbox{div}(\langle \nabla u, X_i(u)\rangle(Y_i(u)-c^i))-\sum_{i=1}^s \langle \partial_t u, X_i(u)\rangle (Y_i(u)-c^i),
\end{eqnarray}
where $c^i\in\mathbb R^{k}$ is an arbitrary constant. Here we have used the killing property of $X_i$ that  yields
$\langle\nabla u, \nabla(X_i(u))\rangle=0$ in the derivation of (\ref{hfhms2}). It is clear that the rest of proof follows exactly
as in Case 1. This completes the proof. \qed\\

\noindent{\bf Proof of Theorem \ref{e-regularity}}. It is easy to see that (\ref{small_norm1}) implies
\begin{equation}\label{small_norm11}
r^{p-(n+2)}\int_{P_r(z)}(|\nabla u|^p+r^p|\partial_t u|^p) \le \epsilon_p^p,
\ \forall z=(x,t)\in P_{\frac{R_0}2}(z_0)\ {\rm{ and }}\ 0<r\le \frac{R_0}2.
\end{equation}
Hence Lemma \ref{grad_estimate5} implies that $u\in C^\infty(P_{\frac{R_0}4}(z_0))$, and
(\ref{gradient_est1}) holds.  \qed

\section{Proof of Theorems \ref{unique1}, \ref{convexity1}, and Corollary \ref{unique_limit}}
\setcounter{equation}{0}
\setcounter{theorem}{0}

In this section, we will provide proofs for our main theorems.  The idea is based on Theorem \ref{e-regularity},
and application of both the Hardy inequality and a generalized Gronwall inequality.\\

\noindent{\bf Proof of Theorem \ref{unique1}}.  For simplicity, we will focus on the case
that $(M,g)$ is a compact Riemannian manifold with boundary and remark on the other two
cases at the end of the proof.

Assume $(M,g)=(\Omega, g_0)$, with $\Omega\subset\mathbb R^n$ and
$g_0$ the standard metric.
By Theorem \ref{e-regularity}, we have that $u_i\in C^\infty(\Omega\times (0, T])$ for $i=1,2$, and
\begin{equation}\label{gradient_est2}
\max\Big\{|\nabla u_1|(x,t), \ |\nabla u_2|(x,t)\Big\}
\le C\epsilon_0\Big(\frac{1}{R_0}+\frac1{d(x,\partial\Omega)}+ \frac{1}{\sqrt{t}}\Big),
\ \forall (x,t)\in \Omega\times (0,T].
\end{equation}
Set $w=u-v$. Then $w$ satisfies
\begin{equation}\label{hfhm1}
\begin{cases}
w_t-\Delta w=A(u)(\nabla u,\nabla u)-A(v)(\nabla v,\nabla v)
&\ {\rm{in}}\ \Omega\times (0,T]\\
w=0 &\ {\rm{on}}\ \partial_p(\Omega\times [0,T]).
\end{cases}
\end{equation}
Multiplying (\ref{hfhm1})
by $w$ and integrating over $\Omega$ yields
\begin{eqnarray*}
\frac{d}{dt}\int_{\Omega}|w|^2+2\int_{\Omega}|\nabla w|^2
&\leq& C\int_{\Omega}(|\nabla u_1|^2+|\nabla u_2|^2)|w|^2
+\int_{\Omega}(|\nabla u_1|+|\nabla u_2|)|\nabla w||w|\\
&\leq& \int_{\Omega}|\nabla w|^2+C\int_{\Omega}(|\nabla u_1|^2+|\nabla u_2|^2)|w|^2.
\end{eqnarray*}
By (\ref{gradient_est2}), the Poincar\'e inequality,  and the Hardy inequality:
$$\int_\Omega \frac{|f(x)|^2}{d^2(x,\partial\Omega)}\lesssim \int_\Omega |\nabla f|^2,
\ \forall f\in H^1_0(\Omega),$$
we have
\begin{eqnarray*}
\frac{d}{dt}\int_{\Omega}|w|^2+\int_{\Omega}|\nabla w|^2
&\le& \frac{C\epsilon_0^2}{R_0^2}\int_\Omega |w|^2+C\epsilon_0^2\int_{\Omega}
\frac{|w(x)|^2}{d^2(x,\partial\Omega)}
+\frac{C\epsilon_0^2}{t}\int_{\Omega}|w|^2\\
&\le& C\Big(\frac{\epsilon_0^2}{R_0^2}+\epsilon_0^2\Big)\int_{\Omega}|\nabla w|^2
+\frac{C\epsilon_0^2}{t}\int_{\Omega}|w|^2.
\end{eqnarray*}
If we choose $\epsilon_0\le (2C)^{-\frac12}$ and $R_0\ge \sqrt{2C}\epsilon_0$, then we have
$C\Big(\frac{\epsilon_0^2}{R_0^2}+\epsilon_0^2\Big)\le 1$ so that
\begin{equation}\label{energy_ineq}
\frac{d}{dt}\int_{\Omega}|w|^2
\le\frac{C\epsilon_0^2}{t}\int_{\Omega}|w|^2.
\end{equation}
This yields
\begin{eqnarray}\label{energy_ineq2}
\frac{d}{dt}\Big(t^{-\frac12}\int_{\Omega}|w|^2\Big)
&=& t^{-\frac12}\frac{d}{dt}\int_{\Omega}|w|^2
-\frac12 t^{-\frac32}\int_{\Omega}|w|^2\nonumber\\
&\le& \left(C\epsilon_0^2-\frac12\right)t^{-\frac32}\int_{\Omega}|w|^2
\le 0.
\end{eqnarray}
Thus we obtain that for any $0<t\le T$,
\begin{equation}\label{gronwall}
t^{-\frac12}\int_{\Omega}|w(x, t)|^2
\le \lim_{s\downarrow 0^+}s^{-\frac12}\int_{\Omega}|w(x, s)|^2.
\end{equation}
Since $w(\cdot,0)=0$, we have
$$w(x,s)=\int_0^s w_\tau(x,\tau)\,d\tau, \ {\rm{a.e.}}\  x\in \Omega$$
so that by the H\"older inequality,
$$s^{-\frac12}\int_{\Omega}|w(x,s)|^2
\le s^\frac12\int_0^s\int_{\Omega}|w_\tau|^2(x,\tau)\,dxd\tau
\le Cs^\frac12\rightarrow 0, \
{\rm{as}}\  s\downarrow 0^+.$$
Thus we conclude that $w\equiv 0$ in $\Omega\times [0,T]$.

When $(M,g)$ is either compact or complete non-compact with $\partial M=\emptyset$, observe that
we can substitute $d(x,\partial M)=\infty$ into the above proof and obtain the same result without
applying the Hardy inequality.  This completes the proof.
\qed\\

\noindent{\bf Proof of Theorem \ref{convexity1}}.  For simplicity, we only consider the difficult case that
$(M,g)$ is compact with boundary. First by Theorem \ref{e-regularity},
we have that $u\in C^\infty(M\times (0, T))$ and
\begin{equation}\label{gradient_est20}
|\nabla u|(x,t)
\le C\epsilon_0\Big(\frac{1}{R_0}+\frac1{d(x,\partial M)}+ \frac{1}{\sqrt{t}}\Big),
\ \forall (x,t)\in M\times (0,T).
\end{equation}

First we need two claims. \\
\noindent{\it Claim 1}. {\it There exists $T_0>0$ such that $\displaystyle\int_{M}|\partial_t u(t)|^2$ is monotone decreasing for $t\ge T_0$:
\begin{equation}
\int_M|\partial_t u|^2(t_2)+C\int_{M\times [t_1, t_2]}|\nabla\partial_t u|^2
\le \int_M|\partial_t u|^2(t_1), \ \ T_0\le t_1\le t_2<T.\label{energy_inequality1}
\end{equation}}
To show it, we introduce the finite quotient for $u$ in the $t$-variable. For sufficiently small $h>0$, set
$$u^h(x,t)=\frac{u(x,t+h)-u(x,t)}{h}, \  \ (x,t)\in M\times (0,T-h).$$
Since  $u^h=0$ on $\partial M$, we see that $u^h\in L^2( [0, T-h], H^1_0(M))$,
$\partial_t u^h \in L^2([0, T-h], L^2(M))$, and
$$\lim_{h\downarrow 0+} \Big\|u^h-u_t\Big\|_{L^2(M\times [0, T-h])}=0.$$
Since $u$ satisfies (\ref{heatflow}), we have
\begin{equation}\label{heatflow10}
u^h_t -\Delta u^h=\frac{1}{h}[A(u(t+h))(\nabla u(t+h), \nabla u(t+h))-A(u(t))(\nabla u(t),\nabla u(t))].
\end{equation}
Multiplying (\ref{heatflow10}) by $u^h$, integrating over $M$, and applying the H\"older inequality and (\ref{gradient_est20}), we obtain
\begin{eqnarray*}
&&\frac12\frac{d}{dt}\int_M |u^h|^2+\int_M |\nabla u^h|^2\\
&\le & C \int_M |u^h|^2(|\nabla u(t+h)|^2+|\nabla u(t)|^2)+|u^h|(|\nabla u(t+h)|+|\nabla u(t)|)|\nabla u^h|\\
&\le & \frac12 \int_M |\nabla u^h|^2 +C\int_M |u^h|^2(|\nabla u(t+h)|^2+|\nabla u(t)|^2)\\
&\le&  \frac12 \int_M |\nabla u^h|^2+C\epsilon_0^2\int_M \Big(\frac{|u^h|^2}{R_0^2}+\frac{|u^h|^2}{d^2(x,\partial M)}+
\frac{|u^h|^2}{t}\Big)\\
&\le&  \frac12 \int_M |\nabla u^h|^2+C\epsilon_0^2\int_M \Big(\frac{|u^h|^2}{R_0^2}+\frac{|u^h|^2}{d^2(x,\partial M)}+
\frac{|u^h|^2}{T_0}\Big)\\
&\le&  \frac12 \int_M |\nabla u^h|^2+C\epsilon_0\int_M |\nabla u^h|^2\le \frac34 \int_M |\nabla u^h|^2,
\end{eqnarray*}
where we have used both the Poincar\'e inequality and the Hardy inequality in the last step, and chosen
$R_0\ge \sqrt{\epsilon_0}$, $T_0\ge \epsilon_0$, and $C\epsilon_0\le \frac14$. Integrating this inequality from $T_0\le t_1\le t_2$ yields
\begin{equation}
\int_M |u^h|^2(t_2)+C\int_{M\times [t_1, t_2]}|\nabla u^h|^2\le \int_M |u^h|^2(t_1).\label{energy_inequality3}
\end{equation}
Sending $h$ to zero in (\ref{energy_inequality3}) yields (\ref{energy_inequality1}).

Next we have\\
\noindent{\it Claim 2}.  {\it There exists $T_0>0$ such that $E(u(t))$ is monotone decreasing for $t\ge T_0$:
\begin{equation}
\int_{M\times [t_1, t_2]}|\partial_t u|^2+E(u(t_2))\le E(u(t_1)), \ \ T_0\le t_1\le t_2<T.
\label{energy_inequality2}
\end{equation} }
For $\delta>0$, let $\phi_\delta\in C_0^\infty(M)$ be a test function such that
$$0\le \phi_\delta\le 1, \ \phi_\delta(x)=1 \ {\rm{for}}\ d(x,\partial M)\ge \delta, \ |\nabla\phi_\delta|\le C\delta^{-1}.$$
Since $u\in C^\infty(M\times (0, T))$, multiplying (\ref{heatflow}) by $\partial_t u\phi_\delta^2$ and
integrating over $M\times [t_1, t_2]$ , we obtain the following local energy inequality:
\begin{equation} \int_{M\times [t_1, t_2]} |\partial_t u|^2\phi_\delta^2+\frac12\int_M |\nabla u(t_2)|^2\phi_\delta^2
\le \frac12\int_M |\nabla u(t_1)|^2\phi_\delta^2+2\int_{M\times [t_1, t_2]} \nabla u \cdot \partial_t u \phi_\delta \nabla\phi_\delta.
\label{energy_inequality4}
\end{equation}
It is clear that (\ref{energy_inequality2}) follows from (\ref{energy_inequality4}), if we can show
\begin{equation}
\label{vanishing_term}
\lim_{\delta\rightarrow 0}\int_{M\times [t_1, t_2]} \nabla u \cdot \partial_t u \phi_\delta \nabla\phi_\delta=0.
\end{equation}
To see (\ref{vanishing_term}), observe that (\ref{energy_inequality1}) implies $\partial_t u(t)\in H^1_0(M)$ for $t\in [t_1,t_2]$
so that
\begin{eqnarray*}
\int_{M\times [t_1, t_2]} |\partial_t  u|^2 |\nabla\phi_\delta|^2
&\lesssim&\delta^{-2}\int_{t_1}^{t_2}\int_{\{x\in M: d(x,\partial M)\le\delta\}}|\partial_t u|^2\\
&\lesssim& \int_{t_1}^{t_2}\int_{\{x\in M: d(x,\partial M)\le\delta\}}|\nabla\partial_t u|^2\rightarrow 0,
\ {\rm{as}}\ \delta\rightarrow 0.
\end{eqnarray*}
It is clear that by the H\"older inequality, (\ref{vanishing_term}) follows from this. Thus (\ref{energy_inequality2}) holds.

Choose $T_0>0$  such that both claims hold. Then by (\ref{heatflow}) we can estimate
\begin{eqnarray}
&& \int_M |\nabla u(t_1)|^2-\int_M |\nabla u(t_2)|^2-\int_M |\nabla (u(t_1)-u(t_2))|^2\nonumber\\
&=&2\int_M \nabla u(t_2)\cdot\nabla (u(t_1)-u(t_2))\nonumber\\
&=& 2\int_M A(u(t_2))(\nabla u(t_2), u(t_2)) \cdot (u(t_1)-u(t_2))
+2\int_M u_t(t_2)\cdot (u(t_2)-u(t_1))\nonumber\\
&=&I+II. \label{convexity30}
\end{eqnarray}
We first estimate $I$.  Recall that for $y\in N$, let $P^\perp(y):\mathbb R^k\to (T_yN)^\perp$ denote the orthogonal
projection from $\mathbb R^k$ to the normal space of $N$ at $y$. Since
$N$ is compact,  a simple geometric argument implies that there exists $C>0$ depending
on $N$ such that
\begin{equation}\label{orth_proj}
\Big|P^\perp(y)(z-y)\Big|\le C|z-y|^2, \ \forall z\in N.
\end{equation}
Thus
\begin{eqnarray*}
|I|&\lesssim& \int_M |\nabla u(t_2)|^2 |u(t_1)-u(t_2)|^2\\
&\leq& C\epsilon_0^2 \int_{M}\Big(\frac{1}{R_0^2}+\frac{1}{T_0}+\frac{1}{d^2(x,\partial M)}\Big)|u(t_1)-u(t_2)|^2\\
&\lesssim&C\epsilon_0\int_{M}|\nabla(u(t_1)-u(t_2))|^2,
\end{eqnarray*}
where we have used both the Poincar\'e inequality and the Hardy inequality in the last step.

By (\ref{energy_inequality1}), we have
$$\int_M |\partial_t u(t_2)|^2\le \frac{1}{t_2-t_1}\int_{M\times [t_1,t_2]}|\partial_t u|^2.$$
This, combined with the H\"older inequality and (\ref{energy_inequality2}), implies
\begin{eqnarray*}
|II|&\leq& \|\partial_t u(t_2)\|_{L^2(M)}\|u(t_1)-u(t_2)\|_{L^2(M)}\\
&\leq & \sqrt{t_2-t_1}\|\partial_t u(t_2)\|_{L^2(M)}\|\partial_t u\|_{L^2(M\times [t_1, t_2])}\\
&\leq &\int_{M\times [t_1,t_2]}|\partial_t u|^2\\
&\le& \frac12 \Big[\int_{M}|\nabla u(t_1)|^2-\int_M |\nabla u(t_2)|^2\Big].
\end{eqnarray*}
Putting the estimates of $I, II$ into (\ref{convexity30}) yields (\ref{convexity2}) so that
the conclusions of Theorem \ref{convexity1} hold. The proof is now complete.
\qed
\\

\noindent{\bf Proof of Corollary \ref{unique_limit}}. It follows from Theorem \ref{convexity1} that
$E(u(t))$ is monotone decreasing for $T_0\le t<+\infty$.
Hence  $$\lim_{t\rightarrow\infty} E(u(t))=c<+\infty.$$
Let $\{t_i\}$ be any monotone increasing sequence such that $\displaystyle\lim_{i\rightarrow\infty} t_i=+\infty$. Then (\ref{convexity2})
implies that for any $j\ge 1$,
$$\int_M |\nabla(u(t_{i+j})-u(t_i))|^2\le C\Big[\int_M |\nabla u(t_i)|^2-\int_M |\nabla u(t_{i+j})|^2\Big]\rightarrow 0$$
as $i\rightarrow\infty$. This implies that there exists a map $u_\infty\in H^1(M, N)$, with $u_\infty=u_0$ on $\partial M$, such that
$$\lim_{t\rightarrow\infty}\int_M|\nabla(u(t)-u_\infty)|^2=0.$$
Since (\ref{energy_inequality2}) implies  there exists $t_i\uparrow\infty$ such that
$$\lim_{i\rightarrow\infty}\|\partial_t u(t_i)\|_{L^2(M)}=0,$$
we see that $u_\infty$ is a weak harmonic map. Moreover, by the gradient estimate (\ref{gradient_est20}), we have that
for any compact set $K\subset\subset M$ and $m\ge 1$, one has that for $t$ sufficiently large,
$$\|\nabla^m u(t)\|_{C^0(K)}\le C(\epsilon_0, m, K),$$
which clearly implies that $u(t)\rightarrow u_\infty$ in $C^m(K)$, as $t\rightarrow\infty$. This completes the proof.
\qed

\section{Serrin's $(l,q)$-solutions and proof of Theorem \ref{unique2}}
\setcounter{equation}{0}
\setcounter{theorem}{0}

In this section, we will indicate that any Serrin's $(l,q)$-solution to (\ref{heatflow}), under a suitable initial-boundary data $u_0$, satisfies the condition (\ref{small_norm1}) for some $p>1$ in Theorem 2.1. We will then sketch a different argument for the
$\epsilon$-regularity, the uniqueness holds for Serrin's $(l,q)$-solution to (\ref{heatflow}) into an arbitrary Riemannian manifold $N$ without boundary.

We start with the following proposition.
\begin{proposition}\label{serrin_prop} For $n\ge 2$, $0<T<+\infty$, and a compact Riemannian manifold
$N\subset\mathbb R^k$ without boundary, suppose $u\in H^1(M\times [0,T],N)$ is a weak solution of (\ref{heatflow}), with
the initial and boundary value $u_0:M\to N$ satisfying $\nabla u_0\in L^r(M)$ for some $n<r<+\infty$,
such that $\nabla u\in L^q_tL^l_x(M\times [0,T])$ for some $(l,q)$ satisfying (\ref{serrin_condition}) with
$l>n, q>2$. Then\\
(i) $\partial_t u\in L^{\frac{q}2}_tL^{\frac{l}2}_x(M\times [0,T])$; and \\
(ii) for any $\epsilon>0$, there exists $R=R(u,\epsilon)>0$ such that
for any $1<s<\min\{\frac{l}2, \ \frac{q}2\}$,
\begin{equation}\label{small_norm6}
\sup\Big\{ r^{s-(n+2)}\int_{P_r(x,t)\cap (M\times [0,T])}(|\nabla u|^s+r^s|\partial_t u|^s)\ | \
(x,t)\in M\times [0, T], \ 0<r\le R\Big\}\le \epsilon^s.
\end{equation}
\end{proposition}
\pf We consider the case that $(M, g)$ is complete and noncompact, and leave the discussion of the other cases to interested readers.
For simplicity, assume $(M, g)=(\mathbb R^n, g_0)$.

Let $H$ be the heat kernel in $\mathbb R^n$. Then by the Duhamel formula, we have
\begin{eqnarray}\label{duhamel}
u(x, t)&=&\int_{\mathbb R^n} H(x-y,t)u_0(y)+\int_0^t\int_{\mathbb R^n} H(x-y, t-s) A(u)(\nabla u,\nabla u)(y,s)\\\
&=& u_1(x,,t)+u_2(x,t).\nonumber
\end{eqnarray}
It is easy to see that
$$\nabla^2 u_1(x,t)=\int_{\mathbb R^n} \nabla_x H(x-y,t)\nabla_y u_0(y).$$
Hence by the standard integral estimates (see \cite{FJR} page 234), we have
\begin{equation}\label{u1-estimate}
\Big\|\nabla^2 u_1\Big\|_{L^{\frac{q}2}_t L^{\frac{l}2}_x(\mathbb R^n\times [0,T])}
\le CT^{\frac12-\frac{n}{2r}}\Big\|\nabla u_0\Big\|_{L^r(\mathbb R^n)}.
\end{equation}
For $u_2$, since
$$\nabla^2 u_2(x,t)=\int_0^t\int_{\mathbb R^n} \nabla_x^2 H(x-y,t-s) A(u)(\nabla u, \nabla u)(y,s),$$
we can apply the Calderon-Zgymund's $L^s_t L^{s'}_x$-theorey to obtain
\begin{equation}\label{u2-estimate}
\Big\|\nabla^2 u_2\Big\|_{L^{\frac{q}2}_t L^{\frac{l}2}_x(\mathbb R^n\times [0,T])}
\le C\Big\||\nabla u|^2\Big\|_{L^{\frac{q}2}_t L^{\frac{l}2}_x(\mathbb R^n\times [0, T])}
\le C\Big\|\nabla u\Big\|_{L^{q}_t L^{l}_x(\mathbb R^n\times [0, T])}^2.
\end{equation}
Substituting (\ref{u1-estimate}) and (\ref{u2-estimate}) into (\ref{duhamel}) yields
$\nabla^2 u \in L^{\frac{q}2}_t L^{\frac{l}2}_x(\mathbb R^n\times [0, T])$. This, combined with the equation
(\ref{heatflow}), then implies (i).

To see (ii), observe that by the H\"older inequality, we have that for any $1<s<\min\{\frac{l}{2}, \frac{q}2\}$,
$$\Big(r^{s-(n+2)}\int_{P_r(x,t)\cap (M\times [0, T])} |\nabla u|^s\Big)^{\frac1{s}}
\le \Big\|\nabla u\Big\|_{L^q_tL^l_x(P_r(x,t)\cap (M\times [0, T]))},$$
and
$$\Big(r^{2s-(n+2)}\int_{P_r(x,t)\cap (M\times [0, T])} |\partial_t u|^s\Big)^{\frac1{s}}
\le \Big\|\partial_t u\Big\|_{L^{\frac{q}2}_tL^{\frac{l}2}_x(P_r(x,t)\cap (M\times [0, T]))}.$$
These two inequalities clearly imply (\ref{small_norm6}).
\qed\\

Now we give a proof of $\epsilon$-regularity of Serrin's solutions to (\ref{heatflow}) for any Riemannian manifold
$N$.
For $x\in\mathbb R^n$, $t>0$,  and $r>0$, let $B_r(x)\subset
\mathbb R^n$ be the ball with center $x$ and radius $r$, and
$P_r(x,t)=B_r(x)\times [t-r^2, t].$  Denote $P_r=P_r(0,0)$.

\blm\label{unq-e-reg-lemma}{\it  There is an $\epsilon_0>0$ such
that if $u\in H^1(P_1,N)$, with $\nabla u\in L^q_tL^l_x(P_1)$
for some $l\geq n$ and $q\geq 2$ satisfying
(\ref{serrin_condition}), is a weak solution to (\ref{heatflow})
and
\beq\label{unq2.1} \|\nabla
u\|_{L^q_tL^l_x(P_1))}\leq\epsilon_0, \eeq
then $u\in C^{\infty}(P_{\frac{1}{2}},N)$
and
\beq\label{unq2.2} \|
u\|_{C^m(P_{\frac{1}{2}})}\leq C(m,n,p,q)\|\nabla
u\|_{L^2(P_1)} \eeq
for any positive integer $m$.}
\elm

We need the following inequality, due to Serrin (\cite{serrin} Lemma 1).

\blm\label{unq-serrin-lemma}{\it For any open set $U\subset\mathbb R^n$ and any open interval $I\subset\mathbb R$,
let $f$, $g$, $h\in
L^2_tH^1_x(U\times I)$ and $f\in L^q_tL^l_x(U\times I)$
with $l\ge n$ and $q\ge 2$ satisfying
(\ref{serrin_condition}). Then we have
\beq\label{unq2.5} \int_{U\times I}|f||g||\nabla h|
\leq
C\|\nabla h\|_{L^2(U\times I)}
\|g\|_{L^2_tH^1_x(U\times I)}^{\frac{n}{l}}
\left\{\int_I\|f\|_{L^l(U)}^{q}\|g\|_{L^2(U)}^2\,dt\right\}^{\frac{1}{q}},
\eeq
where $C>0$ depends only on $n$.
}\elm

\noindent{\bf Proof of Lemma \ref{unq-e-reg-lemma}}.
For any $(x,t)\in P_{\frac{1}{2}}$ and $0<r\leq\frac{1}{2}$,
by (\ref{unq2.1}) we have
\beq\label{unq2.11} \|\nabla
u\|_{L^q_tL^l_x(P_r(x,t)))}\leq\epsilon_0. \eeq
Let $v:
P_r(x,t)\rightarrow \R^k$ solve
 \beq\label{unq2.6} \left\{
\begin{split}
v_t-\Delta v&=0, \quad\mbox{in } P_r(x,t)\\
v&=u, \quad\mbox{on }\partial_pP_r(x,t).
\end{split}\right.
\eeq
Denote $w=u-v$. Multiplying (\ref{heatflow}) and (\ref{unq2.6})
by $w$, subtracting the resulting equations and integrating over
$P_r(x,t)$, we obtain
\beq\label{unq2.7}
\begin{split}
&\sup\limits_{t-r^2\leq s\leq t}\int_{B_r(x)}|w|^2(\cdot,s)+2\int_{P_r(x,t)}|\nabla w|^2
\lesssim \int_{P_r(x,t)}|\nabla u|^2|w|\\
\lesssim &\begin{cases}\|\nabla u\|_{L^2(P_r(x,t))}
\|\nabla w\|_{L^2(P_r(x,t))}^{\frac{n}{l}}
\left\{\int_{t-r^2}^t\|\nabla u\|_{L^l(B_{r}(x))}^{q}\|w\|_{L^2(B_r(x))}^2\right\}^{\frac{1}{q}}, \ & q<\infty\\
\|\nabla u\|_{L^2(P_r(x,t))}
\|\nabla w\|_{L^2(P_r(x,t))}
\|\nabla u\|_{L^\infty L^n(B_{r}(x))},\ & q=\infty
\end{cases}
\end{split}
\eeq
where we have used (\ref{unq2.5}) and the Poincar\'e inequality in last step. Since
$\displaystyle \|\nabla u\|_{L^q_tL^l_x(P_r(z_0))}\le\epsilon$, we obtain,
by the Young inequality,
\beq\label{unq2.8}
\begin{split}
&\sup\limits_{t-r^2\leq s\leq t}\int_{B_r(x)}|w|^2(\cdot,s)+2\int_{P_r(x,t)}|\nabla w|^2\\
\leq&\begin{cases} \|\nabla w\|_{L^2(P_r(x,t))}^{2}+\epsilon_0\|\nabla u\|_{L^2(P_r(x,t))}^2
+C\epsilon_0^{\frac{q}{2}}\|w\|_{L^\infty_tL^2_x(B_r(x))}^2, \ & q<\infty\\
\|\nabla w\|_{L^2(P_r(x,t)}^2+C\epsilon_0^2\|\nabla u\|_{L^2(P_r(x,t))}^2,\ & q=\infty.
\end{cases}
\end{split}\eeq
Choosing $\epsilon_0>0$ so that
$$\begin{cases}
C\epsilon_0^{\frac{q}{2}}\le 1, \ & q<+\infty,\\
C\epsilon_0\le 1, \ & q=\infty,
\end{cases}
$$
we obtain
\beq\label{unq2.10}
\int_{P_r(x,t)}|\nabla w|^2
\leq \epsilon_0\|\nabla u\|_{L^2(P_r(x,t))}^2.
\eeq
On the other hand, by the standard estimate on the heat equation, we obtain that for any $0<\theta<1$,
\beq\label{unq2.20}
(\theta r)^{-n}\int_{P_{\theta r}(x,t)}|\nabla v|^2
\leq C\theta^2r^{-n}\int_{P_r(x,t)}|\nabla u|^2.
\eeq
(\ref{unq2.10}) and (\ref{unq2.20}) imply that
\beq\label{unq2.21}
(\theta r)^{-n}\int_{P_{\theta r}(x,t)}|\nabla u|^2
\leq C\Big(\theta^2+\theta^{-n}\epsilon_0\Big)r^{-n}\int_{P_r(x,t)}|\nabla u|^2.
\eeq
For any $0<\alpha<1$, choose first $\theta_0>0$ such that $C\theta_0^2\leq \frac12\theta_0^{2\alpha}$
 and then
$$\epsilon_0\le \min\left\{\frac{\theta_0^{2\alpha+n}}{2C},
\quad \left(\frac{1}{2C}\right)^{\frac{2}{q}}\right\},$$
we obtain that for any $(x,t)\in P_{\frac{1}{2}}$ and $0<r\leq \frac{1}{2}$, it holds
 \beq\label{unq2.22} (\theta_0
r)^{-n}\int_{P_{\theta_0 r}(x,t)}|\nabla u|^2
\leq\theta_0^{2\alpha}r^{-n}\int_{P_{r_0}(x,t)}|\nabla u|^2.
\eeq
 Iterating (\ref{unq2.22}), we obtain for any positive integer $l$,
 \beq\label{unq2.23} (\theta_0^{l} r)^{-n}\int_{P_{\theta_0^{l}
r}(x,t)}|\nabla u|^2
\leq\theta_0^{2l\alpha}r^{-n}\int_{P_{r}(x,t)}|\nabla
u|^2. \eeq
It is standard that (\ref{unq2.23}) implies
 \beq\label{unq2.24}
r^{-n}\int_{P_{r}(x,t)}|\nabla u|^2
\leq Cr^{2\alpha}\int_{P_1}|\nabla u|^2, \ \forall (x,t)\in P_{\frac12}, \ 0<r\le\frac12.
\eeq
By (\ref{unq2.24}), we have that $\nabla u\in {M}^{2,2-2\alpha}(P_1)$
for any $0<\alpha<1$.
Now we can apply the regularity theorem by  Huang-Wang \cite{huang-wang} Theorem 1.5
to conclude that $u\in C^\infty(P_\frac12)$ and the estimate (\ref{unq2.2}) holds.
This completes the proof.
\endpf\\

By suitable scaling, we have the following estimate on the possible blow-up rate of $\|\nabla
u(t)\|_{L^{\infty}}$ as $t$ tends to zero.

\blm\label{unq-blp-rate-lemma}
For $T>0$ and a compact or complete manifold $(M,g)$
without boundary, suppose that $u\in H^1(M\times[0,T],N)$ is a
weak solution to (\ref{heatflow}), with $\nabla u\in
L^q_tL^l_x(M\times [0,T])$ for some $l>n$ and $q>2$
satisfying (\ref{serrin_condition}), then $u\in
C^{\infty}(M\times(0,T],N)$ and there exists $t_0>0$ such that
 \beq\label{unq2.3}
\sup\limits_{0<t\leq t_0}\sqrt{t}\Big\|\nabla
u(t)\Big\|_{L^{\infty}(M)}\leq C\Big\|\nabla
u\Big\|_{L^q_t L^l_x(M\times [0,t_0])}. \eeq
In particular,
\beq\label{unq2.4}
\lim\limits_{t\downarrow 0^+}\sqrt{t}\Big\|\nabla
u(t)\Big\|_{L^{\infty}(M)}=0. \eeq
\end{lemma}

\pf For simplicity, we assume that $(M,g)=(\mathbb R^n, g_0)$.
Since $\nabla u\in L^q_tL^l_x(\mathbb R^n\times [0,T])$ for some $l>n$ and $q>2$ satisfying (\ref{serrin_condition}),  we have that
for $\epsilon_0>0$ given by Lemma \ref{unq-e-reg-lemma},
there exists $\delta_0>0$ such that
$$ \sup_{(x_0,t_0)\in\R^n\times [0,T]} \Big\|\nabla u\Big\|_{L^q_tL^l_x(P_{\delta_0}(x_0,t_0)\cap \mathbb R^{n+1}_+)}\leq \epsilon_0,$$
In particular, for any $0<\tau\leq \delta_0$ and any $x_0\in\R^n$, we have
\beq\label{unq2.15}
\Big\|\nabla u\Big\|_{L^q_tL^l_x(B_\tau(x_0)\times [0,\tau^2])}\leq \epsilon_0.
\eeq
Define $v(y,s)=u(x_0+\tau y,\tau^2+\tau^2 s)$ for $(y,s)\in P_1(0, 0)$.
Then $v$ solves (\ref{heatflow}) on $P_1(0,0)$, and satisfies
$$
\Big\|\nabla v\Big\|_{L^q_tL^l_x(P_1(0,0))}\leq \epsilon_0.
$$
Hence Lemma \ref{unq-e-reg-lemma} implies
\beq\label{unq2.16} \|\nabla
v\|_{L^{\infty}(P_{\frac{1}{2}}(0,0))}\leq C\|\nabla v\|_{L^2(P_1(0,0))}.
\eeq
After rescalings,  (\ref{unq2.16}) implies that $u\in
C^{\infty}(P_{\frac{\tau}{2}}(x_0,\tau^2))$ and
\beq\label{unq2.17}
\tau\Big\|\nabla u\Big\|_{L^{\infty}(P_{\frac{\tau}{2}}(x,\tau^2))}
\leq C\tau^{-\frac{n}{2}}\Big\|\nabla u\Big\|_{L^2(P_{\tau}(x,\tau^2))}.
\eeq
By H$\ddot{\mbox{o}}$lder's inequality and (\ref{serrin_condition}), we have
\beq\label{unq2.18}
\tau^{-\frac{n}{2}}\Big\|\nabla u\Big\|_{L^2(P_{\tau}(x_0,\tau^2))}
\leq\Big\|\nabla u\Big\|_{L^q_tL^l_x(P_\tau(x_0,\tau^2))}.
\eeq
Putting (\ref{unq2.18}) together with
(\ref{unq2.17}), we obtain
\beq\label{unq2.19}
\tau\Big\|\nabla
u(\tau^2)\Big\|_{L^{\infty}(\mathbb R^n)} \leq C
\Big\|\nabla u\Big\|_{L^q_tL^l_x(\mathbb R^n\times [0,\tau^2])}. \eeq
After sending $\tau\rightarrow 0$, (\ref{unq2.19}) clearly implies (\ref{unq2.4}). It is not hard to see
that (\ref{unq2.3}) also follows.
This completes the proof. \endpf\\

The next lemma handles the case that $(M,g)$ is a compact Riemannian manifold
with boundary.

\begin{lemma}\label{bdry_lemma}
{\it For $T>0$ and a compact manifold $(M,g)$
with boundary, suppose that $u\in H^1(M\times[0,T],N)$ is a
weak solution of (\ref{heatflow}), with $\nabla u\in
L^q_tL^l_x(M\times [0,T])$ for some $l>n$ and $q>2$
satisfying (\ref{serrin_condition}), then $u\in
C^{\infty}(M\times(0,T],N)$. Moreover, for any sufficiently small $\epsilon_0>0$,
there exists $T_0>0$ depending only on $\epsilon_0$ and $u$ such that
\beq\label{unq2.30}
|\nabla u(x_0, t_0)|\leq C\epsilon_0\Big(\frac{1}{d(x_0,\partial M)}+\frac{1}{\sqrt{t_0}}\Big),
\ \forall (x_0,t_0)\in M\times (0,T_0]. \eeq
} \end{lemma}

\pf
Let $\epsilon_0>0$ be given by Lemma \ref{unq-e-reg-lemma}.
Since $\nabla u\in
L^q_tL^l_x(M\times [0,T])$ with $l>n, q>2$, there exists $T_0>0$ such that
$$
\Big\|\nabla u\Big\|_{L^q_tL^l_x(M\times [0,T_0])}
\leq \epsilon_0.$$
For any $x_0\in M$ and $0<t_0\le T_0$, we divide the proof into two cases:\\
(i) $d(x_0,\partial M)>\sqrt{t_0}$; and\\
(ii) $d(x_0,\partial M)\le \sqrt{t_0}$. \\
For (i), since $P_{\sqrt{t_0}}(z_0)\subset M\times (0, T_0]$,
we have $\|\nabla u\|_{L^q_tL^l_x(P_{\sqrt{t_0}}(z_0))}
\le\epsilon_0$. As in Lemma \ref{unq-blp-rate-lemma}, we conclude that
$u\in C^\infty(P_{\frac{\sqrt{t_0}}2}(z_0))$ and
$$|\nabla u|(z_0)\le \frac{C\epsilon_0}{\sqrt{t_0}}.$$
For (ii), set $r_0=\min\{d(x_0,\partial M), \sqrt{t_0}\}$. Then
$P_{r_0}(z_0)\subset M\times (0, T_0]$ and $\|\nabla u\|_{L^q_tL^l_x(P_{r_0}(z_0))}
\le\epsilon_0$. Hence we can conclude that
$u\in C^\infty(P_{\frac{r_0}2}(z_0))$ and
$$|\nabla u|(z_0)\le \frac{C\epsilon_0}{r_0}\le C\epsilon_0\Big(\frac{1}{d(x_0,\partial M)}+\frac{1}{\sqrt{t_0}}\Big).$$
Thus (\ref{unq2.30}) holds. This completes the proof. \qed \\

\noindent{\bf Proof of Theorem \ref{unique2}}. It follows from Lemma \ref{unq-blp-rate-lemma} and Lemma \ref{bdry_lemma} that
there exists $T_0>0$ such that both the condition (\ref{struwe_mono}) of Theorem \ref{unique1}
and the estimate (\ref{gradient_est1}) of Theorem \ref{e-regularity} hold on $M\times [0,T_0]$. Thus we can apply the same proof
of Theorem \ref{unique1} to obtain that $u=v$ on $M\times [0, T_0]$.
One can repeat the same argument to show that $u=v$ on $M\times [T_0, T]$.\qed

\section{Appendix}

\setcounter{equation}{0}
\setcounter{theorem}{0}

As a byproduct of the proof of Theorem \ref{unique1}, we will prove a convexity property
on certain weak harmonic maps  that yields an alternative, simple proof of
the  uniqueness theorem on the Dirichlet problem of weak harmonic maps,  due
to Struwe \cite{struwe2} for $n=3$ and Moser \cite{moser1} for $n\ge 4$. Furthermore,
the statement of the uniqueness theorem for $N$ either a unit sphere or a compact Riemannian homogeneous
manifold without boundary is an improvement of that by \cite{struwe2} and \cite{moser1}.

To do it, we introduce the Morrey spaces in $\mathbb R^n$.  For $1\le l<+\infty$, $0<\lambda\le n$,
$0<R\le +\infty$, and $U\subset \mathbb R^n$, $f\in  M^{l,\lambda}_R(U)$ iff $f\in L^l_{\hbox{loc}}(U)$
satisfies
$$
\|f\|_{M^{l,\lambda}_R(U)}^l
:=\sup_{x\in U}\sup_{0<r\le \min\{R, d(x,\partial U)\}}
\Big\{r^{\lambda-n}\int_{B_r(x)}|f|^l\Big\}<+\infty.$$
Denote $M^{p,\lambda}(U)=M^{p,\lambda}_\infty(U)$.

For any bounded smooth domain $\Omega\subset \mathbb R^n$, we have
\begin{theorem}\label{convexity_hm} For $n\ge 2$, $\delta\in (0,1)$, and $1<p\le 2$,
there exist $\epsilon_p=\epsilon(p,\delta)>0$ and $R_p=R(p,\delta)>0$ such that
if $u\in H^{1}(\Omega, N)$ is a weak
harmonic map satisfying either\\
(i) $\displaystyle\|\nabla u\|_{M^{2,2}_{R_2}(\Omega)} \le\epsilon_2$, when $N$ is a
compact Riemannian manifold without boundary, or \\
(ii) $\displaystyle\|\nabla u\|_{M^{p,p}_{R_p}(\Omega)} \le\epsilon_p$,
when $N=S^{k-1}$ or a compact Riemannian homogeneous manifold without boundary.
Then
\begin{equation}\label{convex_ineq}
\int_{\Omega}|\nabla v|^2\ge \int_{\Omega}|\nabla u|^2+(1-\delta)\int_{\Omega}|\nabla(v-u)|^2
\end{equation}
holds for any $v\in H^1(\Omega, N)$ with $v=u$ on $\partial\Omega$.
\end{theorem}
\pf
First, as observed by \cite{struwe2} and \cite{moser1},   for an arbitrary manifold $N$ under the condition (i),
the small energy regularity theorem on stationary harmonic maps by Bethuel
\cite{bethuel} holds. While, for $N=S^{k-1}$ under the condition (ii), the small energy regularity theorem on
weak harmonic maps by Moser \cite{moser2} or Lemma 2.3 is applicable.
Thus  we have $u\in C^\infty(\Omega, N)$ and, for any $x\in\Omega$,  it holds
\begin{equation}\label{grad_estimate}
|\nabla u|(x)
\le {C\epsilon_p}(\frac{1}{d(x,\partial\Omega)}+\frac{1}{R_p}).
\end{equation}
Here $p=2$ for an arbitrary $N$.

Now multiplying the equation of $u$ by ($u-v$) and integrating over $\Omega$, we obtain
\begin{equation}\label{test_eqn}
\int_{\Omega}\nabla u\cdot\nabla (u-v)
=\int_{\Omega}\langle A(u)(\nabla u,\nabla u), u-v\rangle.
\end{equation}
This, combined with (\ref{grad_estimate}), the Poincar\'e inequality, and the Hardy inequality,  implies
\begin{eqnarray}
\Big|\int_{\Omega}\langle A(u)(\nabla u,\nabla u), u-v\rangle\Big|
&\le& C\int_{\Omega}|\nabla u|^2|u-v|^2\nonumber\\
&\le& C\epsilon_p^2\int_{\Omega}\frac{|u-v|^2}{R_p^2}+\frac{|u-v|^2}{d(x,\partial\Omega)^2}\nonumber\\
&\le& C\epsilon_p^2(1+\frac{1}{R_p^2})\int_{\Omega}|\nabla(u-v)|^2\nonumber\\
&\le & \frac{\delta}{2} \int_{\Omega}|\nabla(u-v)|^2\label{small_rhs}
\end{eqnarray}
provided that we have chosen $\displaystyle\epsilon_p\le (\frac{\delta}{4C})^\frac12$ and $R_p$ such
$\displaystyle\frac{C\epsilon_p^2}{R_p^2}\le\frac{\delta}4$.
Thus, by (\ref{test_eqn}) and (\ref{small_rhs}) we obtain
\begin{eqnarray*}
&&\int_{\Omega}|\nabla v|^2-\int_{\Omega}|\nabla u|^2-\int_{\Omega}|\nabla(v-u)|^2\\
&&=2\int_{\Omega}\nabla u\cdot\nabla(v-u)=-2\int_{\Omega}\langle A(u)(\nabla u,\nabla u), u-v\rangle\\
&&\ge -\delta\int_{\Omega}|\nabla(v-u)|^2.
\end{eqnarray*}
This clearly implies (\ref{convex_ineq}), provided that $\epsilon>0$ is sufficiently small. This proof is complete. \qed

\begin{corollary}\label{uniqueness} For
$n\ge 2$ and $1<p\le 2$,  there exist $\epsilon_p>0$
and $R_p>0$ such that if $u_1, u_2\in H^{1}(\Omega, N)$ are two weak harmonic maps satisfying either\\
(i) $\displaystyle\max_{i=1}^2\|\nabla u_i\|_{M^{2,2}_{R_2}(\Omega)}\le\epsilon_2$, when $N$ is a
compact Riemannian manifold without boundary, or \\
(ii) $\displaystyle\max_{i=1}^2\|\nabla u_i\|_{M^{p,p}_{R_p}(\Omega)}
\le\epsilon_p$,
when $N=S^{k-1}$ or a compact Riemannian homogeneous manifold without boundary.\\
Then $u_1\equiv u_2$ in $\Omega$,  provided that $u_1-u_2\in W^{1,2}_0(\Omega,\mathbb R^k)$.
\end{corollary}
\pf  Choosing $\delta=\frac12$, we can apply Theorem \ref{convexity_hm} to $u_1$ and $u_2$
by choosing sufficiently small $\epsilon_p>0$ and $R_p>0$. Thus Theorem  \ref{convexity_hm}
implies
$$
\int_{\Omega}|\nabla u_2|^2\ge \int_{\Omega}|\nabla u_1|^2+\frac12\int_{\Omega}|\nabla(u_2-u_1)|^2,
$$
and
$$
\int_{\Omega}|\nabla u_1|^2\ge \int_{\Omega}|\nabla u_2|^2+\frac12\int_{\Omega}|\nabla(u_1-u_2)|^2.
$$
Adding these two inequalities together yields
$$\int_{\Omega}|\nabla(u_1-u_2)|^2=0.$$
Therefore, $u_1\equiv u_2$ in $\Omega$.\qed

\noindent{\bf Acknowledgements}.  Both authors are partially supported by NSF grant 1001115. The second author is also
partially supported by NSFC grant 11128102.


\begin{thebibliography}{99}

\bibitem{bethuel}
F. Bethuel, {\em On the singular set of stationary harmonic maps}.
Manuscripta Math., {\bf 78} (1993), 417-443.

\bibitem{BBC}
F. Bethuel,  J. Coron, J.  Ghidaglia, A. Soyeur,
{\em Heat flows and relaxed energies for harmonic maps}.
Nonlinear diffusion equations and their equilibrium states, 3 (Gregynog, 1989), 99-109, Progr. Nonlinear Differential Equations Appl., 7, Birkhäuser Boston, Boston, MA, 1992.

\bibitem{cam} S. Campanato,
{\em Equazioni ellittiche del $II^0$ ordine espazi $\mathcal L^{(2,\lambda)}$.}
Ann. Mat. Pura Appl., (4) {\bf 69} (1965) 321-381.

\bibitem{chang}

 K. Chang, {\em Heat flow and boundary value problem for harmonic maps}.
Ann. Inst. H. Poincar\'e Anal. Non Lin\'eaire, {\bf 6} (5)
(1989), 363-395.

\bibitem{chen-ding}

Y. Chen, W. Ding, {\em Blow-up and global existence for heat flows
of harmonic maps}. Invent. Math., {\bf 99} (3) (1990), 567-578.

\bibitem{chang-ding-ye}

 K. Chang, W. Ding, R. Ye, {\em Finite-time blow-up of the heat flow of
harmonic maps from surfaces}.  J. Diff. Geom., {\bf 36} (1992),
507-515.

\bibitem{CWY}
A. Chang,  L. Wang, P. Yang,
{\em Regularity of harmonic maps.}
Comm. Pure Appl. Math., {\bf 52} (1999), no. 9, 1099-1111.

\bibitem{chen-lin}

Y. Chen, F. Lin, {\em Evolution of harmonic maps with Dirichlet
boundary conditions}. Comm. Anal. Geom., {\bf 1}(3-4) (1993), 327-346.

\bibitem{CLL}

Y. Chen, J. Li, F. Lin, {\em Partial regularity for weak heat flows into spheres}.
Comm. Pure Appl. Math., {\bf 48} (1995), no. 4, 429-448.


\bibitem{chen-struwe}

Y. Chen, M. Struwe, {\em Existence and partial regularity results
for the heat flow for harmonic maps}. Math. Z., {\bf 201} (1) (1989), 83-103.

\bibitem{CW}

Y. Chen, C. Wang, {\em Partial regularity for weak heat flows into Riemannian homogeneous spaces}.
Comm. Partial Differential Equations, {\bf 21} (1996), no. 5-6, 735-761.

\bibitem{coron}

J. Coron, {\em Nonuniqueness for the heat flow of harmonic maps}.
Ann. Inst. H. Poincar\'e Anal. Non Lin\'eaire, {\bf 7} (1990), no. 4, 335-344.



\bibitem{ES}
J. Eells, J. Sampson, {\em Harmonic mappings of Riemannian manifolds}.
Amer. J. Math., {\bf 86} (1964), 109-160.

\bibitem{evans}
L. C. Evans, {\em Partial regularity for stationary harmonic maps into spheres.} Arch. Rational Mech. Anal.,
{\bf 116} (1991), no. 2, 101-113.

\bibitem{feldman}

M. Feldman, {\em Partial regularity for harmonic maps of evolution into spheres}.
Comm. Partial Differential Equations, {\bf 19} (1994), no. 5-6, 761-790.

\bibitem{FJR}

E. Fabes, F. Jones, N. Riviere, {\em  The initial value problem for the Navier-Stokes equations
with date in $L^p$}. Arch. Rational Mech. Anal., {\bf 45} (1972), 222-240.

\bibitem{freire}

A. Freire, {\em Uniqueness for the harmonic map flow from surfaces
to general targets}. Comm. Math. Helv., {\bf 70} (1) (1995), 310-338.



\bibitem{H}
R. Hamilton, {Harmonic maps of manifolds with boundary.}
Lecture Notes in Mathematics, Vol. 471. Springer-Verlag, Berlin-New York, 1975.

\bibitem{helein0}
F. H\'elein, {\em R\'egularit\'e des applications faiblement harmoniques entre une surface et une sph\'ere.}
(French) [Regularity of weakly harmonic maps between a surface and an n-sphere]
C. R. Acad. Sci. Paris S\'er. I Math., {\bf 311} (1990), no. 9, 519-524.

\bibitem{helein} F. H\'elein, {\em Regularity of weakly harmonic maps from a surface into a manifold with symmetries.}
Manuscripta Math., {\bf 70} (1991), no. 2, 203-218.

\bibitem{HHW} J. Hineman, T. Huang, C. Y. Wang, {\em Regularity and  uniqueness of heat flow of biharmonic maps.} Preprint, 2012.


\bibitem{HKW}

S. Hildebrandt, H. Kaul, K. Widman, {\em An existence theorem for harmonic mappings of Riemannian manifolds}.
Acta Math., {\bf 138} (1977), no. 1-2, 1-16.

\bibitem{hartman} P. Hartman,
{\em On homotopic harmonic maps.}
Canad. J. Math., {\bf 19} (1967) 673-687.

\bibitem{huang-wang}

T. Huang, C. Y. Wang,
{\em Notes on the regularity of harmonic map systems}.
Proc. Amer. Math. Soc., {\bf 138} (6) (2010), 2015-2023.

\bibitem{JN} F. John, L. Niernberg,
{\em On functions of bounded mean oscillation.} Comm. Pure Appl. Math.,
{\bf 14} (1961), 415-426.



\bibitem{lin-wang}

F. H. Lin, C. Y. Wang, {\em On the uniqueness of heat flow of
harmonic maps and hydrodynamic flow of nematic liquid crystals}.
Chin. Ann. Math., {\bf 31B} (6) (2010), 921-938.


\bibitem{lin} L. Z. Lin,  {\em Uniformity of harmonic map heat flow at infinite time}. Preprint (2012),  arXiv:1202.5756.




\bibitem{moser}

R. Moser, {\em Regularity for the approximated harmonic map equation and application to the heat flow for harmonic maps}.
Math. Z., {\bf 243} (2003), no. 2, 263-289.

\bibitem{moser1}

R. Moser, {\em Unique solvability of the Dirichlet problem for weakly harmonic maps}.
Manuscripta Math., {\bf 105} (2001), no. 3, 379-399.


\bibitem{moser2}
R. Moser, {\em An $\epsilon$-regularity result for generalized harmonic maps into spheres}. Electron. J. Differential Equations, (2003), No. 1, 7 pp.

\bibitem{moser3} R. Moser, {\em An $L^p$ regularity theory for harmonic maps}. Preprint (2012).


\bibitem{RS} T. Rivier\`e, M. Struwe,
{\em Partial regularity for harmonic maps and related problems.} Comm. Pure Appl. Math.,
{\bf 61} (2008), no. 4, 451-463.

\bibitem{schoen}
R. Schoen, {\em Analytic aspects of the harmonic map problem.}
Seminar on nonlinear partial differential equations (Berkeley, Calif., 1983), 321-358, Math. Sci. Res. Inst. Publ., 2. Springer, New York-Berlin, 1984.

\bibitem{serrin}

J. Serrin, {\em The initial value problem for the Navier-Stokes
equations}.  In Nonlinear Problems (Proc. Sympos., Madison, Wis.,
pages 69-98. Univ. of Wisconsin Press, Madison, Wis., 1963).

\bibitem{simon} L. Simon,
{\em Asymptotics for a class of nonlinear evolution equations, with applications to geometric problems.}
Ann. of Math., (2) {\bf 118} (1983), no. 3, 525-571.

\bibitem{struwe}

M. Struwe, {\em On the evolution of harmonic mappings of Riemannian
surfaces}. Comm. Math. Helv., {\bf 60} (1985), 558-581.

\bibitem{struwe1}
M. Struwe, {\em On the evolution of harmonic maps in higher dimensions}.
J. Diff. Geom., {\bf 28} (1988), no. 3, 485-502.

\bibitem{struwe2}
M. Struwe, {\em Uniqueness of harmonic maps with small energy}.
Manuscripta Math., {\bf 96} (1998), no. 4, 463-486.

\bibitem{struwe3}
M. Struwe, {\em Geometric evolution problems}.
Nonlinear partial differential equations in differential geometry (Park City, UT, 1992), 257-339,
IAS/Park City Math. Ser., 2, Amer. Math. Soc., Providence, RI, 1996.

\bibitem{topping}
P. Topping, {\em Rigidity in the harmonic map heat flow.}
J. Differential Geom., {\bf 45} (1997), no. 3, 593–610.

\bibitem{wang}

 C. Y. Wang,
{\em Heat flow of harmonic maps whose gradients belong to $L^n_xL^\infty_t$}.
Arch. Rational Mech. Anal., {\bf 188} (2008), 309-349.

\bibitem{wang-well}

C. Y. Wang, {\em Well-posedness for the heat flow of harmonic maps
and the liquid crystal flow with rough initial data}.  Arch. Rational
Mech. Anal., {\bf 200} (2011), 1-19.

\bibitem{luwang}

L. Wang, {\em Harmonic map heat flow with rough boundary data}.
Trans. Amer. Math. Soc., {\bf 364} (2012), 5265-5283.

\end{thebibliography}
\end{document}